\documentclass[11pt]{article}
\usepackage{graphics,latexsym,amssymb,amsmath,amscd}
\RequirePackage{graphics}
\usepackage{graphicx}
\usepackage{psfrag}
\input epsf
\def\abs#1{\left \vert #1 \right \vert} 
\def\Im{\hbox{\rm Im}\,}
\def\[#1\]{\begin{eqnarray*}#1\end{eqnarray*}}
\def\Re{\hbox{\rm Re}\,}

\def\pu{{\bf PU}(2,1)}

\def\ch#1{{{\bf H}^{#1}_{\C}}}
\def\rh#1{{{\bf H}^{#1}_{\R}}}
\def\Heis{{\mathfrak N}}
\def\isom{{\widehat{\bf PU}(2,1)}}
\def\phi{\varphi}

\newtheorem{thm}{Theorem}[section]
\newtheorem{dfn}[thm]{Definition}

\newtheorem{cor}[thm]{Corollary}
\newtheorem{prop}[thm]{Proposition}
\newtheorem{lem}[thm]{Lemma}

\newcommand{\Pf}{{\sc Proof}. }
\newcommand{\EPf}{\hbox{}\hfill$\Box$\vspace{.5cm}}

\def\A{{\mathbb A}}

\input xy
\xyoption{all}

\newcommand{\C}{{{\mathbb C}}}
\newcommand{\R}{{{\mathbb R}}}
\newcommand{\Z}{{{\mathbb Z}}}
\newcommand{\PP}{{{\mathbb P}}}
\begin{document}
\title{A spherical CR structure on the complement
of the figure eight knot with discrete holonomy.
\small{}}
\author{ E. Falbel\\
Institut de Math{\'e}matiques\\
Universit{\'e} Pierre et Marie Curie\\
4, place Jussieu \\ F-75252 Paris\\ 
e-mail: {\tt falbel\@math.jussieu.fr}}
\maketitle
\begin{abstract}
We describe a general geometrical construction of spherical CR structures.
We construct then spherical CR structures
on the complement of the figure eight knot and the Whitehead link.  They 
have discrete holonomies contained
in $PU(2,1,\Z[\omega])$ and $PU(2,1,\Z[i])$ respectively.  These are the same 
ring of integers appearing in the real hyperbolic geometry of the
 corresponding links.
\end{abstract}

\section{Introduction}

One of the most important  examples of hyperbolic
manifolds is the complement of the figure eight knot.
It was shown by Riley in \cite{R} that the fundamental group of that
manifold had a discrete representation in $PSL(2,\C)$.  In fact 
he showed that there exists a representation contained
in $PSL(2,\Z[\omega])$ where  $\Z[\omega]$ is the ring 
of Eisenstein integers.
On the other hand the construction by Thurston is based 
on gluing of ideal tetrahedra and that led to general constructions
on a large family of 3-manifolds.

It is not known which hyperbolic manifolds 
admit a spherical CR structure.  In fact very few constructions
of 
spherical CR 3-manifolds with discrete holonomy
 exist at all.  The only construction of such
 a structure on a  3-manifold (which is not a circle
bundle) previous to this work
  is essentially 
for the Whitehead link and other manifolds obtained from it by 
Dehn surgery 
in \cite{S1,S2}.

We propose a geometrical construction by gluing appropriate
tetrahedra adapted to CR geometry. {\bf In particular we prove
in this paper that the complement of the figure eight knot 
has a  spherical CR structure 
with discrete holonomy such that the holonomy of the boundary torus is 
 parabolic and faithful} (see Theorem \ref{main} and Proposition \ref{torus}).
   As another example we  also construct a 
spherical structure on the complement of the Whitehead link
which differs from \cite{S1} with discrete holonomy (Theorem \ref{main1}).
  It is interesting to observe that
we obtain representations of the fundamental groups of those link complements
with  values in $PU(2,1,\Z[\omega])$ and 
$PU(2,1,\Z[i])$, that is the same rings of integers of the complete
structures in the case of real hyperbolic geometry.

There are two different aspects in the construction.  The first
is a very general method to construct CR manifolds by gluing
which can be used to construct CR structures on many other
(hyperbolic or not) manifolds.  The second aspect is discreteness
which in the real  hyperbolic case is much simpler to decide than in
the CR case because of the absence of a metric structure in the latter. 
  For the complement
of the figure eight knot and the Whitehead link we show discreteness
of the representation of the fundamental group by explicitly showing
that the group is a subgroup of $PU(2,1,\Z[\omega])$ and $PU(2,1,\Z[i])$
respectively.

We thank R. Benedetti, M. Deraux, W. Goldman, J.-P. Koseleff, J. Parker, J. Paupert, R. Schwartz and P. Will for many fruitful
discussions.

\section{Complex hyperbolic space}

\subsection{$\pu$, $\isom$ and the Heisenberg group}\label{stereographic} 

Let $\C^{2,1}$ denote the complex vector space equipped with the Hermitian
form 
$$
\langle z,w\rangle = z_1\overline{w}_3+z_2\overline{w}_2+z_3\overline{w}_1.
$$
Consider the following subspaces in $\C^{2,1}$:
\[
V_+ &=& \{ z\in {\C}^{2,1}\ \ :\ \ \langle z,z\rangle > 0 \ \},\\
V_0 &=& \{ z\in {\C}^{2,1}\setminus\{ 0\}\ \ 
:\ \ \langle z,z\rangle = 0 \ \},\\
V_- &=& \{ z\in {\C}^{2,1}\ \ :\ \ \langle z,z\rangle < 0 \ \}.
\]
Let $P:{\C}^{2,1}\setminus\{ 0\} \rightarrow {\C}P^{2}$ be the
canonical projection onto complex projective space. Then
$\ch{2} = P(V_-)$ equipped with the Bergman metric is complex hyperbolic 
space. The boundary of complex hyperbolic space is $P(V_0) = \partial\ch{2}$. 
The isometry group $\isom$ of $\ch{2}$ comprises holomorphic 
transformations in $\pu$, the unitary group of $\langle \cdot,\cdot\rangle$, 
and anti-holomorphic transformations arising elements of $\pu$ followed by
complex conjugation.

\def\H{{\bf H}}
The {\sl Heisenberg group} $\Heis$ is the set of pairs $(z, t)\in
{\C}\times{\R}$ with the product 
$$
(z,t)\cdot (z',t') = (z+z', t + t' + 2 \Im z \overline{z}').
$$
Using stereographic projection, we can identify $\partial\ch{2}$
with the one-point compactification $\overline{\Heis}$ of $\Heis$. 
The Heisenberg group acts on itself by left translations. Heisenberg 
translations by $(0, t)$ for $t\in\R$ are called {\sl vertical translations}.

Define the inversion in the $x$-axis in $\C\subset\Heis$ by
$$
\iota_x:(z,t) \mapsto (\overline {z},-t).
$$
All these actions extend trivially to the compactification
$\overline{\Heis}$ of $\Heis$ and represent transformations in $\isom$
acting on the boundary of complex hyperbolic space (see \cite{G}).

A point $p=(z,t)$ in the Heisenberg group and the point $\infty$ are lifted 
to the following points in $\C^{2,1}$:
$$
\hat{p}=\left[\begin{matrix}\frac{-|z|^2+it}{2} \\ z \\ 1 \end{matrix}\right]
\quad\hbox{ and }\quad
\hat{\infty}=\left[\begin{matrix} 1 \\ 0 \\ 0 \end{matrix}\right].
$$
Given any three points $p_1$, $p_2$, $p_3$ in $\partial\ch{2}$ we define 
{\sl Cartan's angular invariant} $\A$ as 
$$
\A(p_1,p_2,p_3)=
\arg(-\langle\hat{p}_1,\hat{p}_2\rangle\langle\hat{p}_2,\hat{p}_3\rangle
\langle\hat{p}_3,\hat{p}_1\rangle).
$$
In the special case where $p_1=\infty$, $p_2=(0,0)$ and $p_3=(z,t)$ we simply 
get $\tan(\A)=t/\abs{z}^2$. 

\subsection{$\R$-circles, $\C$-circles and $\C$-surfaces }

There are two kinds of totally geodesic submanifolds of real dimension 2
in $\ch{2}$: complex lines in $\ch{2}$ are {\sl complex geodesics} 
(represented by $\ch{1} \subset \ch{2}$) and Lagrangian planes in
$\ch{2}$ are {\sl totally real geodesic 2-planes}
(represented by $\rh{2}\subset \ch{2}$). Each of these totally
geodesic submanifolds is a model of the hyperbolic plane. 

%A discrete subgroup of $\pu$ preserving a complex line is called
%{\sl $\C$-Fuchsian} and is isomorphic to a subgroup of
%${\bf P}\bigl({\bf U}(1)\times{\bf U}(1,1)\bigr)\subset \pu$.
%A discrete subgroup of $\pu$ preserving a Lagrangian plane is called
%{\sl $\R$-Fuchsian} and is isomorphic to a subgroup of
%${\bf SO}(2,1)$ included in $\pu$ by the projectivisation of the
%obvious inclusion ${\bf SO}(2,1)\subset{\bf SU}(2,1)$.

Consider complex hyperbolic space $\ch{2}$ and its boundary 
$\partial \ch{2}$. We define {\sl {$\C$}-circles} in $\partial\ch{2}$ to 
be the boundaries of complex geodesics in $\ch{2}$.
Analogously, we define {\sl {$\R$}-circles} in $\partial\ch{2}$ to be
the boundaries of Lagrangian planes in $H_{\C}^2$.

\begin{prop}[see \cite{G}]
In the Heisenberg model, $\C$-circles are either vertical lines or
ellipses, whose projection on the $z$-plane are circles.
\end{prop}

Finite $\C$-circles are determined by a {\sl centre} $M=(z=a+ib ,c)$ and a
{\sl radius} $R$. They may also be described using polar vectors in $P(V_+)$
(see Goldman \cite{G} page 129).

 If we use the Hermitian form 
$\langle \cdot,\cdot\rangle$, a finite chain with centre $(a+ib,c)$ and 
radius $R$ has polar vector (that is the orthogonal vector in $\C^{2,1}$ to
the plane determined by the chain).
$$
\left[\begin{matrix} \frac{R^2-a^2-b^2+ic}{2} \\ a+ib \\ 1 \end{matrix}\right]
$$

Given two points $p_1$ and $p_2$ in Heisenberg space, we write
$[p_1,p_2]$ for a choice of one of the two segments of $\C$-circle
joining them.  The choice will be determined from the context.

\begin{dfn}
A $\C$-triangle determined by three points
$[p_0,p_1,p_2]$ is a triangular surface determined by 
segments of $\C$-circles joining $p_0$ to each point
a segment of $\C$-circle $[p_1,p_2]$.
\end{dfn}
Observe that, in principle, there are four smooth triangular
surfaces canonically associated to $[p_0,p_1,p_2]$.  Each of those triangles
could be part of a $\C$-sphere (see \cite{FZ}).

\section{Tetrahedra}

To copy the tetrahedra of the conformal case we start with 4 points such that
each triple of points up to a sign has a fixed  Cartan's invariant.
The edges of the tedrahedron could be segments of either $\R$-circles
 or $\C$-circles and the faces 
should be adapted later to that one skeleton.  In this paper we will use
$\C$-circles and $\C$-triangles.
%\begin{figure}[]
%\center{{\scalebox{0.5}{\includegraphics{tetra.eps}}}}
%\caption{{}}
%\label{}
%\end{figure}
%\begin{figure}[]
%\center{{\scalebox{0.5}{\includegraphics{Rtetra.eps}}}}
%\caption{{}}
%\label{}
%\end{figure}

\subsection{CR triples of points}
We first describe triples of points in the standard spherical CR sphere.
They are classified up to $PU(2,1)$ in the following proposition.
\begin{prop}[\cite{C}, see \cite{G}]
The Cartan invariant classifies triples of points up to $PU(2,1)$.
\end{prop}

A natural way to obtain a triangle is then to join the 3 points by 
$\C$-circles. As in spherical geometry, for each pair of points there 
are two choices of circular segments joining them.

\subsubsection{CR tetrahedra}
For a general tetrahedra we have 4 Cartan invariants corresponding 
to each triple
of points.  But one of them is determined by the others in view 
of the cocycle condition (see \cite{G} pg. 219):
$$
-A(x_2,x_3,x_4)+A(x_1,x_3,x_4)-A(x_1,x_2,x_4)+A(x_1,x_2,x_3)=0.
$$
As a special case of tetrahedra we have the following.
\begin{prop}
If three triples of four points are contained in $\R$-circles ($\C$-circles),
 the  four points are contained in a common $\R$-circle ($\C$-circle).
\end{prop}
\Pf We will prove the result on $\R$-circles, the other case being easier.
 From the cocycle relation, each triple is contained in 
an $\R$-circle as $A=0$ for all triples.  Without loss of generality,
we can suppose that three of the points are $\infty, [0,0], [1,0]$
in Heisenberg coordinates.  The fourth point is in an $\R$-circle
containing $\infty, [0,0]$ on one hand, so it is in the plane $t=0$.
 On the other hand, it should be in an $\R$-circle passing through
$[1,0]$ and $\infty$, that is in the contact plane at $[1,0]$. 
 The intersection of both planes is precisely the $x$-axis.
\EPf

\begin{dfn}A tetrahedron is a configuration of four points and a choice of
edges, that is a choice of $\C$-circle segments joining
each pair of points.
\end{dfn}

\begin{dfn} A symmetric tetrahedron
 is a configuration of four points with an anti-holomorphic
symmetry and
a choice of $\C$-circle segments joining
each pair of points.
\end{dfn} 
By normalizing the coordinates of the four points
we can assume that they are given by
$$  p_1=\infty\ \  p_2=0\ \  q_1=(1,t)\ \  q_2=(z,s|z|^2)
$$

\begin{lem}[cf. \cite{W}]The configuration of four points $p_1$,$p_2$
,$q_1$ and $q_2$ has a $\Z_2$ anti-holomorphic symmetry exchanging 
$p_1, p_2$ and $q_1, q_2$ if and
only if $t=s$.
\end{lem}
\Pf A simple proof  follows writing
the general form of an anti-holomorphic transformation permuting
$\infty$ and $0$.  It is given by 
$$(z,t)\rightarrow \left ( -\frac{\bar z}{|\lambda|^2(|z|^2+it)},
\frac{t}{|\lambda|^4(|z|^4+t^2)}\right ),
$$where $\lambda\in \C^*$.  Imposing that the points $q_1$ and $q_2$
are permuted then gives the result.

\EPf

In that case $\A(p_1,p_2,q_1)=\A(p_1,p_2,q_2)$ and  
$\A(p_1,q_1,q_2)=\A(p_2,q_1,q_2)$. 
\begin{lem} For configurations with $\Z_2$ symmetry as above,
 $\A(p_1,p_2,q_1)=\A(p_1,q_1,q_2)$ if and only if $tg(\A(p_1,p_2,q_1))= t=\frac{\Im z}{1-\Re z}$.
\end{lem}
\Pf A simple computation shows that
$$\A(p_1,p_2,q_1)=arg(|z|^2(1+it)/2)=arctg (t)$$ and 
$$\A(p_1,q_1,q_2)=arg(|z|^2(1+it)/2 +(1-2\bar{z}-it)/2)
=arctg\frac{t(|z|^2-1)+2\Im z}{|z|^2+1-2\Re z}.$$
The proposition follows by equating the two formulas and solving for $t$.

\EPf

\begin{dfn}
We call a symmetric tetrahedron regular if the configuration of four points
satisfies
$$
\A(p_1,p_2,q_1)=\A(p_1,p_2,q_4)=\A(p_1,q_1,q_2)=\A(p_2,q_1,q_2).
$$
\end{dfn}
In that case $s=t$ and $t=\frac{\Im z}{1-\Re z}$.

%\subsection{Special Ideal tetrahedra}

%A special symmetric tetrahedron is defined
% by the points $p_1=(0,n)$, $p_2=(0,-n)$, $q_1=(\omega,0)$
%and $q_2=(1,0)$, with $\omega=e^{i\theta}$.
%We will compute the angular invariants in that case.
%  Observe that there is a symmetry
%which exchanges $p_1$ and  $p_2$ at the same time as  $q_1$ and $q_2$.
%\begin{lem}
%The $\R$-reflection on the infinite $\R$ circle passing through $[0,0]$
%and $[ e^{i\theta/2},0]$ is a symmetry of the tetrahedron which exchanges $p_1$ and  $p_2$ at the same time as  $q_1$ and $q_2$.
%\end{lem} 
% %
%This shows that special tetrahedra have a $\Z_2$ symmetry;
% they are symmetric and moreover 
% the $\C$-circles determined by 
%$[p_1,p_2]$ and $[q_1,q_2]$ are orthogonal.

%The complex inversion 
%$$
%(z,t)\rightarrow \left ( \frac{z}{|z|^2-it},-\frac{t}{|z|^4+t^2}\right)
%$$
%fixes the equator and interchanges the poles $p_1$ and $p_2$.
%\begin{lem}
%\begin{itemize}
%\item  $A(q_1,q_2,p_2)=-A(q_2,q_1,p_1)=\arg (1-\omega)$
%\item  $A(p_1,q_2,p_2)= -A(p_2,q_1,p_1)=arctg \frac{1-n^2}{2n}$
%\end{itemize}
%\end{lem}
%\Pf $A(q_1,q_2,p_2)= \arg(-\langle\hat{q}_1,\hat{q}_2\rangle\langle
%\hat{q}_2,\hat{p}_2\rangle
%\langle\hat{q}_2,\hat{q}_1\rangle )=\arg ( -(-2+2\omega)(in-1)(-in-1))=
%\arg (1-\omega)$.  On the other hand
% $A(p_1,q_2,p_2)= \arg(-\langle\hat{p}_1,\hat{q}_2\rangle\langle
%\hat{q}_2,\hat{p}_2\rangle
%\langle\hat{q}_2,\hat{p}_1\rangle )=\\
%\arg ( -(in-1)(-1+in)(-in-in))=arctg \frac{1-n^2}{2n}$.  
%\EPf

\subsubsection{Parameters of ideal tetrahedra}

\begin{figure}\label{Figure:regular}
\setlength{\unitlength}{1cm}
\begin{picture}(13,10)
%\put(3.6,4.5){$q_1$}\put(4.4,5.5){$\tilde z_2$}\put(4.8,4.7){$\tilde z_1$}\put(4.4,3.7){$\tilde z_3$}
%\put(8.3,.8){$p_2$}\put(8.3,2){$z_1$}\put(9.4,2){$z_2$}\put(6.4,1.9){$z_3$}
%\put(12.1,4.5){$q_2$}\put(11.5,5.5){$\tilde z_3$}\put(10.8,4.7){$\tilde z_1$}\put(11.5,3.7){$\tilde z_2$}
%\put(8.4,9){$p_1$}\put(8.3,8){$z_1$}\put(9.4,8){$z_3$}\put(6.5,8){$z_2$}
\put(4,1){\includegraphics[height=8cm,width=8cm]{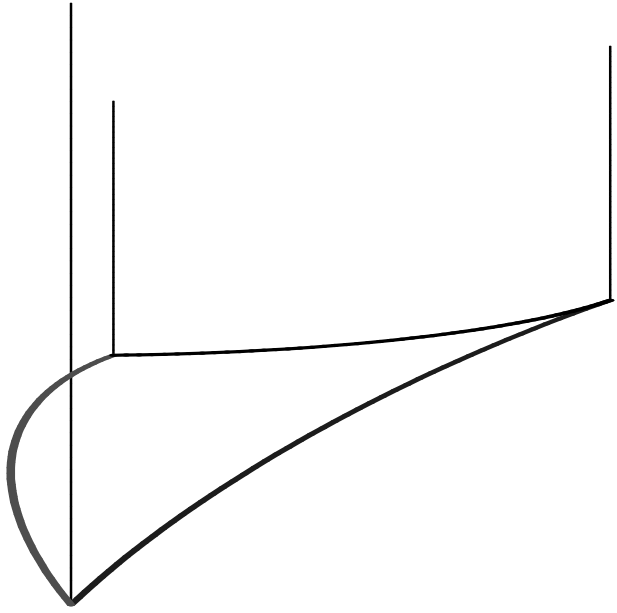}}
\end{picture}
\caption{\sl The standard tetrahedron}
\end{figure}

To each vertex of a tetrahedra we associate the complex
 coordinates of the three
vertical lines obtained when we place that vertex at $\infty$.  That gives
us four euclidean triangles. 
Consider the following configuration of points
$$  p_1=\infty\ \  p_2=0\ \  q_1=(1,t)\ \  q_2=(z,s|z|^2).
$$
There are many possible choices for edges.  We will choose the infinite
edges to be the halves of the  vertical $\C$-circles which tend to $+\infty$.
The other edges will be clear from the context.  In particular
for $$  p_1=\infty\ \  p_2=0\ \  q_1=(1,\sqrt{3})\ \  q_2=(\frac{1}{2}+\frac{i\sqrt{3}}{2},\sqrt{3}),
$$ we have Figure \ref{Figure:regular}.
The invariant of the triangle determined by the points
$(p_2,q_1,q_2)$ at the line determined by
$p_2$ is $z$.
In order to obtain the invariant of the triangle determined by the triple 
$(q_2,p_1,p_2)$ at $q_2$ 
we use the complex
inversion 
$$
I(z,t)=(\frac{z}{|z|^2-it},\frac{-t}{|z|^4+t^2})
$$ to move the point $q_1$ to $\infty$.  We proceed in the same manner 
for the other points.

\begin{figure}
\setlength{\unitlength}{1cm}
\begin{center}
\begin{picture}(7,10)
\psfrag{z1}{$z_1$}
\psfrag{z1'}{$z_1'$}
\psfrag{z~1}{$\tilde z_1$}
\psfrag{z~1'}{$\tilde z_1'$}
\psfrag{z2}{$z_2$}
\psfrag{z2'}{$z_2'$}
\psfrag{z~2}{$\tilde z_2$}
\psfrag{z~2'}{$\tilde z_2'$}
\psfrag{z3}{$z_3$}
\psfrag{z3'}{$z_3'$}
\psfrag{z~3}{$\tilde z_3$}
\psfrag{z~3'}{$\tilde z_3'$}
\psfrag{q1}{$q_1$}
\psfrag{q2}{$q_2$}
\psfrag{p1}{$p_1$}
\psfrag{p2}{$p_2$}
{\scalebox{.8}{\includegraphics[height=8cm,width=8cm]{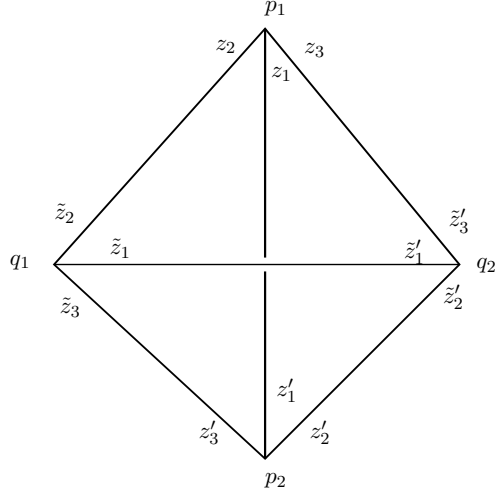}}}
\end{picture}
\end{center}
\caption{\sl Parameters for a CR tetrahedron}\label{Figure:parameter}
\end{figure}

We consider Figure \ref{Figure:parameter} to describe the
parameters of a tetrahedron. Note that, contrary to the ideal
tetrahedron in real hyperbolic geometry, the euclidean  invariant 
at each vertex is not the same. The following proposition follows
immediately from the considerations above by a simple calculation.
\begin{prop}
 For a tetrahedron  given by 
$$  p_1=\infty\ \  p_2=0\ \  q_1=(1,t)\ \  q_2=(z,s|z|^2)
$$then $z_1=z$,$z_1'=\frac{i+t}{\bar{z}(i+s)}$,
$\tilde{z}_1=z\frac{t+i-\bar{z}(i+s)}{(z-1)(t-i)}$ and $\tilde{z}_1'=
\frac{1}{\bar{{{z}}}}\frac{-(i+t) +\bar{z}(i+s)}{(z-1)(i-s)}$.
Where, as usual, $z_2=\frac{1}{1-z_1}$ and  $z_3=1-\frac{1}{z_1}$ and so on.
$$
tg \A(p_1,p_2,q_1)= t
$$
$$
tg \A(p_1,q_1,q_2)= \frac{|z_1|^2s-t+2\Im{z_1}}{|z_1-1|^2}
$$
$$
tg\A(p_1,p_2,q_2)=s
$$
\end{prop}
Here, the three Cartan invariants are independent.
We also have the following relations
$$
\tilde{z}=\frac{z(z'-1)(t+i)}{z'(z-1)(t-i)}\ \ {\mbox \rm and}\ \  
\tilde{z}'=\frac{(z'-1)(i+s)}{(z-1)(i-s)}.
$$
Therefore
$$
t=i\frac{zz'-z-\tilde{z}z'+\tilde{z}z'z}{-zz'+z-\tilde{z}z'+\tilde{z}z'z}
\ \ {\mbox \rm and}\ \  
s=i\frac{z'-1-\tilde{z}'+\tilde{z}'z}{-z'+1-\tilde{z}'+\tilde{z}'z}.
$$

For the symmetric
tetrahedra, the situation is simpler: 
\begin{cor}
For a symmetric tetrahedron  given by 
$$  p_1=\infty\ \  p_2=0\ \  q_1=(1,t)\ \  q_2=(z,t|z|^2)
$$then $z_1=z$,$z_1'=z/|z|^2$,
$\tilde{z}_1=z\frac{(\bar{z}-1)(1-it)}{({z}-1)(1+it)}$ and $\tilde{z}_1'=\tilde{z}_1/|\tilde{z}_1|^2$.
Where, as usual, $z_2=\frac{1}{1-z_1}$ and  $z_3=1-\frac{1}{z_1}$.
$$
tg A(p_1,p_2,q_1)= -i\frac{ z_1\bar z_1-z_1 -z_1\tilde z_1+ \tilde z_1}
{ z_1\bar z_1-z_1 +z_1\tilde z_1-\tilde z_1}=t
$$
$$
tg A(p_1,q_1,q_2)= -i\frac{ z_1\bar z_1+z_1 -z_1\tilde z_1 -\tilde z_1}
{ z_1\bar z_1-z_1 +z_1\tilde z_1-\tilde z_1}=
\frac{t(|z|^2-1)+2\Im z}{|z-1|^2}
$$
\end{cor}
The proposition above shows that the set of symmetric tetrahedra is
parametrized by a strictly pseudoconvex CR hypersurface in $\C\times\C$,
 namely, 
solving for $t$, we obtain the equation

%$$
%it_3= \frac{ -(z_1-1)\tilde z_1+ z_1(\bar z_1-1)}{(z_1-1)\tilde z_1+ z_1(\bar z_1-1)}
%$$
%and therefore, for $(z_1-1)\tilde z_1+ z_1(\bar z_1-1)\neq 0$ 
%(otherwise it corresponds to degenerate tetrahedra),
$$
|z_1|=|\tilde z_1|.
$$

%One sees that this is a strictly pseudoconvex CR  hypersurface.  Special
%tetrahedra are precisely those with module $|z_1|=|\tilde z_1 |=1$ (degenerate
%tetrahedra correspond to $\tilde z_1=1$):
\begin{prop}
If the special symmetric tetrahedron is given by 
$$  p_1=(0,t)\ \  p_2=(0,-t)\ \  q_1=(1,0)\ \  q_2=(e^{i\theta},0)
$$then $z_1=e^{i\theta}$ and $\tilde{z}_1=\frac{(t+i)^2}{(t-i)^2}$.
Where, as usual, $z_2=\frac{1}{1-z_1}$ and  $z_3=1-\frac{1}{z_1}$.
\end{prop}

In the regular symmetric case there is only one complex parameter
which should be compared to the  parameter for real hyperbolic
tetrahedra:
\begin{prop}
Regular symmetric tetrahedra are parametrized by the
complex number $z_1=\tilde z_1=z$ with $\Re z\neq 1$.
In the coordinates above, $tg (\A(p_1,p_2,q_1))=tg( A(p_1,q_1,q_2))= t=\frac{\Im z}{1-\Re z}$.
\end{prop}

This proposition shows that if a real hyperbolic ideal triangulation
has modular invariants for its tetrahedra contained in a line
$\frac{\Im z}{1-\Re z}=constant$ gives rise to representations
of the fundamental group of the the manifold into PU(2,1).

As a last observation,
 the  moduli for a tetrahedron can be expressed using other invariants
as the Koranny-Reimann cross-ratio and Cartan's invariant.
\subsubsection{The standard special tetrahedron}
  We make $\omega=e^{-i\pi/3}$ and
$t=2+\sqrt{3}$ for a special tetrahedra.  Using the formulas above we obtain
\begin{lem} If $p_1=(0,2+\sqrt{3})$, $p_2=(0,-(2+\sqrt{3}))$, $q_1=(\omega,0)$ 
and $q_2=(1,0)$ then in the parameters above $z_1=\tilde z_1=\bar\omega$.  
Moreover the tetrahedron is symmetric and 
$A(q_1,q_2,p_2)=\frac{\pi}{3}$ and
  $A(p_1,q_2,p_2)= -\frac{\pi}{3}$.
\end{lem}

\subsubsection{Another special tetrahedron}\label{sub:whitehead}

 We make $p_1=(0,1+\sqrt{2})$, $p_2=(0,-(1+\sqrt{2}))$, $q_1=(1,0)$
and $q_2=(i,0)$.  We obtain the following

\begin{lem}For $T_w=[p_1,p_2,q_1,q_2]$ as above $z_1=\tilde z_1=i$.
\end{lem}

\subsection{Fundamental lemma for special symmetric tetrahedra}
 We define the procedure of filling the faces from the one skeleton
of the  tetrahedra in such a way that the
 2-skeleton will be $\Z_2$-invariant:

\begin{dfn}
The diverging $\C$-rays procedure is the definition of the 2-skeleton
 by taking $\C$-segments from $p_1$ to the edges $[q_1,q_2]$, $[q_2,p_2]$
and $\C$-segments from $p_2$ to the edges $[q_1,q_2]$, $[q_1,p_1]$.
\end{dfn}
Observe that the rays start from $p_1$ or $p_2$ and not from $q_1$ or $q_2$.
 
\begin{lem} The special symmetric
 tetrahedron defined by the procedure of diverging 
$\C$-rays is homeomorphic to a tetrahedron.
\end{lem}\label{lemma:fundamental}
\Pf We make one of the vertexes go to infinity keeping the other
on the vertical axis.  The other pair of points is in an orthogonal 
$\C$-circle.  They correspond to the normalization
$p_1=\infty$, $p_2=0$, $q_1=(1,t_3)$ and $q_2=(e^{i\theta},t_3)$.
The computations are easy and show that the faces don't
intersect.

\EPf

\section{Gluing the standard tetrahedron: figure eight knot}
\begin{figure}
\setlength{\unitlength}{1cm}
\begin{picture}(13,8)
\put(4.8,4.5){$q_1$}
\put(8.3,.8){$p_2$}
\put(12.1,4.5){$q_3$}
\put(8.4,9){$p_1$}
\put(10.5,7.2){$g_1$}
\put(4,3){$g_3$}
\put(8.8,3.6){$g_2$}
\put(4,1){\includegraphics[height=8cm,width=8cm]{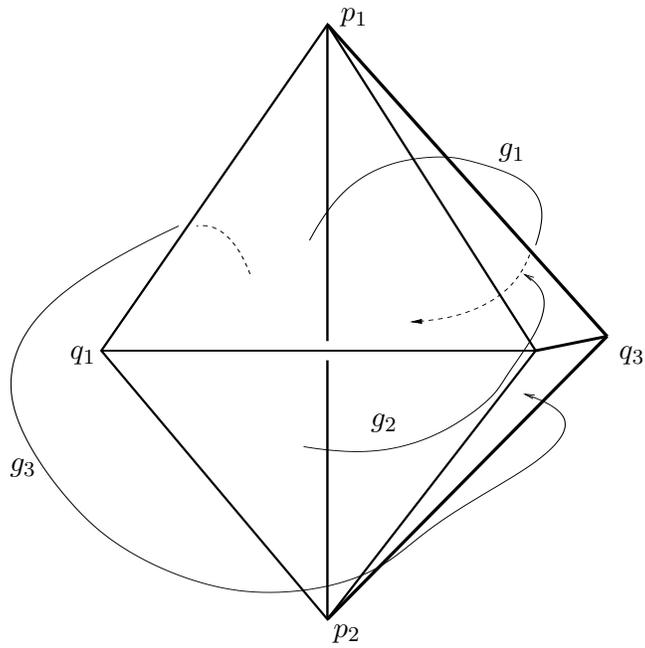}}
\end{picture}
\caption{\sl Identifications on the tetrahedra.}\label{figure:tetrahedron2}
\end{figure}

\begin{thm}
There exists a  spherical $CR$-structure on the complement of the 
figure eight knot with discrete holonomy.
\end{thm}
\Pf We use the same identifications that Thurston used 
in his construction for a hyperbolic real structure on the figure
eight knot.  That is, two tetrahedra with the  identifications
given in Figure \ref{figure:tetrahedron2}. We realize the two tetrahedra in
the Heisenberg space gluing a pair of sides.
The side pairings transformations are shown in Figure \ref{figure:tetrahedron2}
where the two tetrahedra are represented with a common side (here we introduce the point
 $q_3=(\bar\omega,0)$).  They are determined by their action on three points and are
defined by:
$$
g_1 :(q_2,q_1,p_1)\rightarrow (q_3,p_2,p_1)
$$
$$
g_2 :(p_2,q_1,q_2)\rightarrow (p_1,q_3,q_2)
$$
$$
g_3 :(q_1,p_2,p_1)\rightarrow (q_2,p_2,q_3)
$$
%%%%%%%%%%%%%%%%%%%%%%%%%%%%%%%%%%%%%%%%%%%%%%%%%%%%%%%%%%%%%%%%%
%\begin{figure}[]
%\center{{\scalebox{.5}{\includegraphics{eight.eps}}}}
%\caption{Identifications on two tetrahedra to obtain the figure eight knot}
%\label{figure:eight}
%\end{figure}
%%%%%%%%%%%%%%%%%%%%%%%%%%%%%%%%%%%%%%%%%%%%%%%%%%%%%%%%%%%%%%%%%

\begin{figure}
\setlength{\unitlength}{1cm}
\begin{picture}(13,10)
\put(1.8,4.5){$q_1=(\omega,0)$}
\put(8.3,.8){$p_2=(0,-2-\sqrt{3})$}
\put(12.1,4.5){$q_2=(1,0)$}
\put(8.4,9){$p_1=(0,2+\sqrt{3})$}
\put(4,1){\includegraphics[height=8cm,width=8cm]{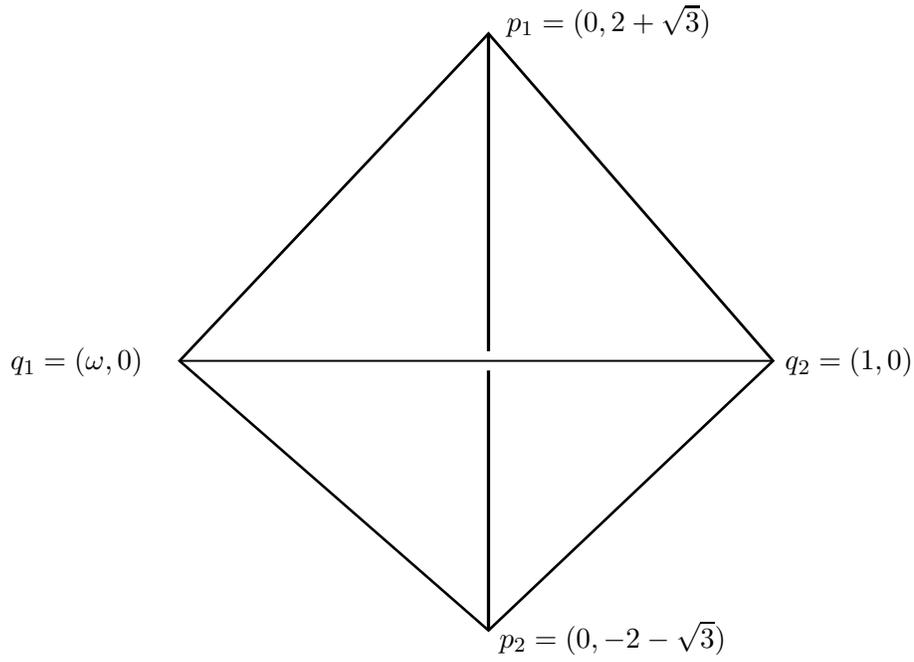}}
\end{picture}
\caption{\sl A schematic view of the standard ideal
tetrahedron in the Heisenberg group}
\end{figure}
%%%%%%%%%%%%%%%%%%%%%%%%%%%%%%%%%%%%%%%%%%%%%%%%%%%%%%%%%%%%%%%%%%%%
%\begin{figure}[]
%\center{{\scalebox{.5}{\includegraphics{tetrahedron.eps}}}}
%\caption{A schematical view of the ideal tetrahedron.}
%\label{figure:tetrahedron}
%\end{figure}

\begin{figure}[]
\center{{\scalebox{.7}{\includegraphics{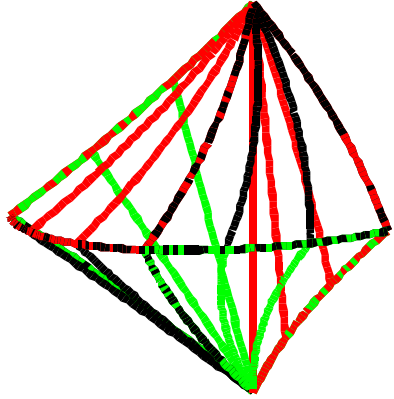}}}}\center{{\scalebox{.7}{\includegraphics{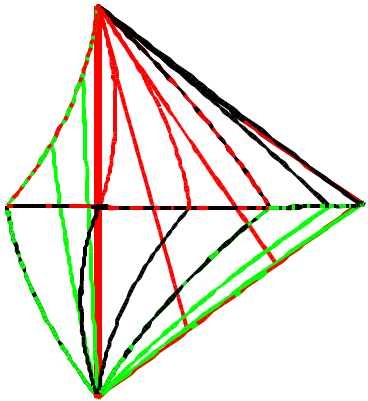}}}}
\caption{Identification of two tetrahedra to obtain the figure eight knot in coordinates of the Heisenberg group: two views.}
\label{figure:eightheiseinberg}
\end{figure}

%\begin{figure}[]
%\center{{\scalebox{.5}{\includegraphics{tetrahedron2.eps}}}}
%\caption{Identifications on the tetrahedra.}
%\label{figure:tetrahedron2}
%\end{figure}
In order to define 
in a compatible way the faces we join the vertex $p_1$ to each point in 
the edge
$[q_1,q_2]$ with segments of $\C$-circles 
and use $g_1$ to define the corresponding face. In the same
manner, we define the other two pairs of identified faces.
We verify that the faces are compatible and that
 the structure is well defined around the edges.  This follows 
because the triangles at the vertexes are equilateral. 
  The drawings in Figure \ref{figure:eightheiseinberg}
show the faces.  

The rest of the proof concerns information about the holonomy
 of the structure, we divide it
in several subsections. 
\subsection{Discreteness of the representation}

Recall from above and  Figure \ref{figure:tetrahedron2} the side pairing transformations 
of the two tetrahedron with a common side:
$$
g_1 :(q_2,q_1,p_1)\rightarrow (q_3,p_2,p_1)
$$
$$
g_2 :(p_2,q_1,q_2)\rightarrow (p_1,q_3,q_2)
$$
$$
g_3 :(q_1,p_2,p_1)\rightarrow (q_2,p_2,q_3)
$$
We conjugate each generator by the map
$$
\gamma : (\infty,0,[1,-\sqrt{3}])\rightarrow (p_1,q_2,q_1)
$$
and obtain after some computation
the following matrices in $SU(2,1,\Z[\omega])$, the Eisenstein-Picard group (see \cite{FP}):
$$
G_1=\left[\begin{matrix} 1  &  \omega  &  -\omega\\
                         0  &   1      &  -\bar\omega\\
		 	 0  &   0 &  1
\end{matrix}\right]
$$
$$
G_2=\left[\begin{matrix} 1  &  1  &  -\omega\\
                         -1  &   0      &    -\bar\omega\\
		 	-\bar\omega   &   \omega &  1
\end{matrix}\right]
$$
$$
G_3=\left[\begin{matrix} 1  &  1  &  -\omega\\
                         -\omega  &   \bar\omega      & -1 -\bar\omega\\
		 	-\bar\omega   &   0 &  1+\omega
\end{matrix}\right]
$$
Note that $G_1$, $G_3$ are parabolic and $G_2$ is elliptic.
\begin{thm}\label{main}
The fundamental group of the complement of the figure eight knot
 has a discrete representation in $PU(2,1)$.
\end{thm}
\Pf As the generators are in  $SU(2,1,\Z[\omega])$, the group is discrete.
\EPf

\subsection{Holonomy of the torus link}

\begin{figure}
\setlength{\unitlength}{1cm}
\begin{picture}(13,10)
\psfrag{Id}{Id}
\psfrag{a}{a}
\psfrag{b}{b}
\psfrag{c}{c}
\psfrag{d}{d}
\psfrag{e}{e}
\psfrag{f}{f}
\psfrag{g}{g}
\psfrag{h}{h}
\psfrag{G_1}{$G_1$}
\psfrag{G_2}{$G_2$}
\psfrag{G_3}{$G_3$}
\put(3,1){\includegraphics[height=8cm,width=12cm]{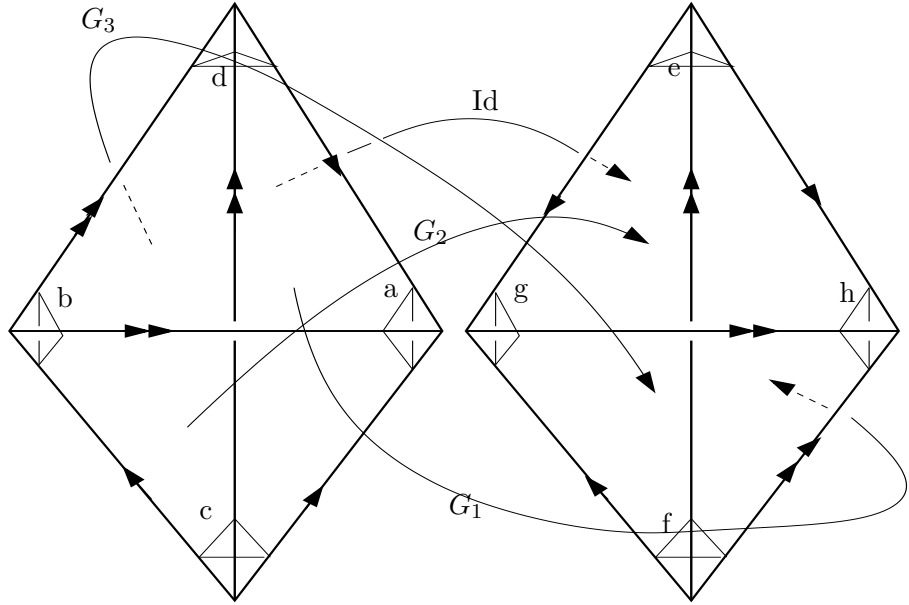}}
\end{picture}
\caption{\sl The triangulation of the torus link}\label{figure:link}
\end{figure}

\begin{figure}
\setlength{\unitlength}{1cm}
\begin{picture}(13,10)
\psfrag{Id}{Id}
\psfrag{a}{a}
\psfrag{b}{b}
\psfrag{c}{c}
\psfrag{d}{d}
\psfrag{e}{e}
\psfrag{f}{f}
\psfrag{g}{g}
\psfrag{h}{h}
\psfrag{g_1}{$G_1$}
\psfrag{g_2}{$G_2$}
\psfrag{g_3}{$G_3$}
\psfrag{g_2^{-1}}{$G_2^{-1}$}
\psfrag{g_1^{-1}}{$G_1^{-1}$}
\psfrag{g_3^{-1}}{$G_3^{-1}$}
\psfrag{tz1}{$\tilde z_1$}
\psfrag{tz3}{$\tilde z_3$}
\psfrag{tz2}{$\tilde z_2$}
\psfrag{w1}{$w_1$}
\psfrag{w2}{$w_2$}
\psfrag{w3}{$w_3$}
\psfrag{z1}{$z_1$}
\psfrag{z3}{$z_3$}
\psfrag{z2}{$z_2$}
\psfrag{tw1}{$\tilde w_1$}
\psfrag{tw3}{$\tilde w_3$}
\psfrag{tw2}{$\tilde w_2$}
\put(1,1){\includegraphics[height=8cm,width=16cm]{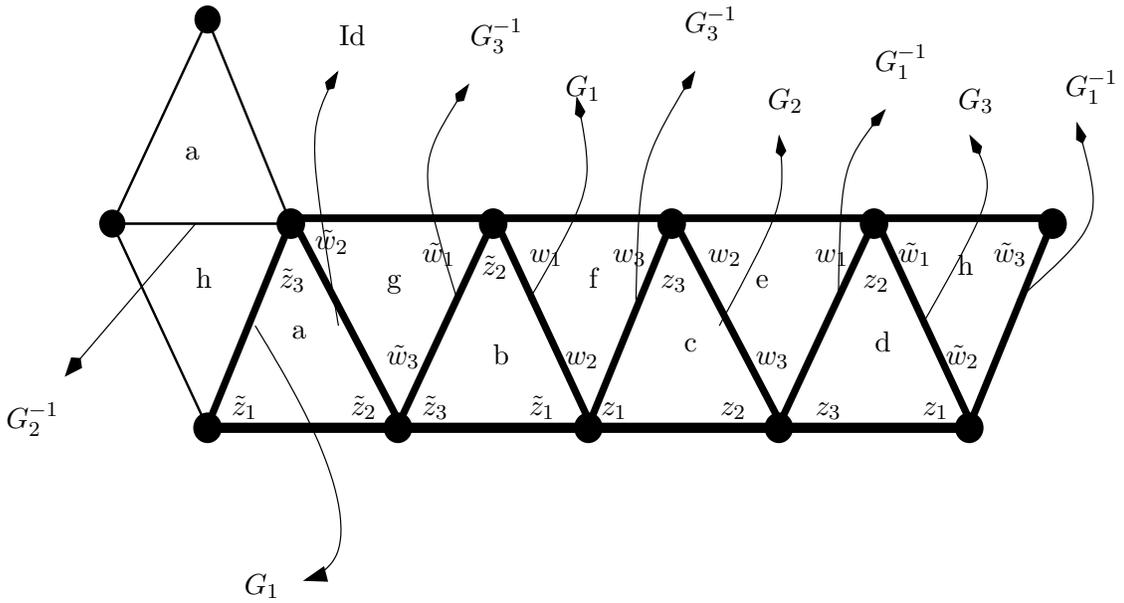}}
\end{picture}
\caption{\sl Computation of the holonomy at the vertex}\label{figure:holo}
\end{figure}

Refering to Figure \ref{figure:link},
the holonomy of the torus link at the vertex can be computed following
the identifications of the triangles forming the link. Starting with the
triangle on the right of the first tetrahedron we obtain the generators
$$
H_1=G_1^{-1}G_3G_1^{-1}G_2G_3^{-1}G_1G_3^{-1}=\left[\begin{matrix} 1  &  0  &  0\\
                        -2\bar\omega  &   1      &  0\\
		 	-2\omega-1  &   2\omega &  1
\end{matrix}\right]
$$
$$
H_2=G_2^{-1}G_1=\left[\begin{matrix} 1  &  0  &  0\\
                         \bar\omega  &   1      &  0\\
		 	 -\omega  &   -\omega &  1
\end{matrix}\right]
$$
Figures  \ref{figure:link} and
\ref{figure:holo} shows how to compute those elements. 
 It turns out
that they are parabolic and independent:
\begin{prop}\label{torus}
The holonomy of the torus link is faithful and parabolic.
\end{prop}
Moreover $H_1$ is a Heisenberg translation by $[2\omega,2\sqrt{3}]$ and
$H_2$ is a Heisenberg translation by $[-\omega,\sqrt{3}]$.  Therefore
the holonomy is generated by 
$[-\omega,\sqrt{3}]$ ($H_2$) and a vertical translation
 $[0,4\sqrt{3}]$ ($H_1 H_2^2$).  Of course,
that should be the case as the group is commutative and discrete.

\subsection{Relation to Eisenstein-Picard group}
In \cite{FP} we proved that the Eisenstein-Picard Group $PU(2,1,\Z[\omega])$
 is generated by
$$
P=\left[\begin{matrix}
1 & 1      & -\omega  \\
0 & -\omega & +\omega \\
0 & 0      & 1
\end{matrix} \right],\quad
Q=\left[ \begin{matrix}
1 &  1 & -\omega \\
0 & -1 & 1 \\
0 & 0  & 1
\end{matrix} \right].
$$
and
$$
I=\left[ \begin{matrix}
0 &  0 & 1 \\
0 & -1 & 0 \\
1 & 0  & 0
\end{matrix} \right].
$$
In this section we identify the generators of the holonomy in terms 
of these generators.  The information is contained in the
following proposition.  We first state a lemma whose
 proof is a simple computation after
a guess obtained by identifying the translational part
of each parabolic element. 
\begin{lem}
The holonomy of the torus link is given by
$$
H_1 = I(QP^{-1}Q(PQ^{-1})^{-2})^2I
$$
and
$$
H_2 = IPQ^{-1}P^2Q^{-1}I
$$
\end{lem}
From the lemma and a computation we obtain the generators
of the group.
\begin{prop}
$G_1= PQ^{-1}P^2Q^{-1}$,
$G_2= I (PQ^{-1}P^2Q^{-1})^{-1}I PQ^{-1}P^2Q^{-1}$ and
$G_3= A I H_2 I A^{-1}$ where
$A= PQ^{-2}(PQ^{-1})^2I (QP^{-1})^2 P$
\end{prop}

\section{Equations along the edges}

\begin{figure}
\setlength{\unitlength}{1cm}
\begin{picture}(13,10)
\psfrag{Id}{Id}
\psfrag{w=q_1}{$w=q_1$}
\psfrag{w=q_3}{$w=q_3$}
\psfrag{v=q_4}{$v=q_4$}
\psfrag{v=q_2}{$v=q_2$}
\psfrag{v=p_1}{$v=p_1$}
\psfrag{v=p_2}{$v=p_2$}
\psfrag{w}{w}
\psfrag{a}{a}
\psfrag{b}{b}
\psfrag{c}{c}
\psfrag{d}{d}
\psfrag{e}{e}
\psfrag{f}{f}
\psfrag{g}{g}
\psfrag{h}{h}
\psfrag{g_1}{$G_1$}
\psfrag{g_2}{$G_2$}
\psfrag{g_3}{$G_3$}
\psfrag{g_2^{-1}}{$G_2^{-1}$}
\psfrag{g_1^{-1}}{$G_1^{-1}$}
\psfrag{g_3^{-1}}{$G_3^{-1}$}
\psfrag{tz1}{$\tilde z_1$}
\psfrag{tz3}{$\tilde z_3$}
\psfrag{tz2}{$\tilde z_2$}
\psfrag{w1}{$w_1$}
\psfrag{w2}{$w_2$}
\psfrag{w3}{$w_3$}
\psfrag{z1}{$z_1$}
\psfrag{z3}{$z_3$}
\psfrag{z2}{$z_2$}
\psfrag{tw1}{$\tilde w_1$}
\psfrag{tw3}{$\tilde w_3$}
\psfrag{tw2}{$\tilde w_2$}
\put(3,1){\includegraphics[height=5cm,width=10cm]{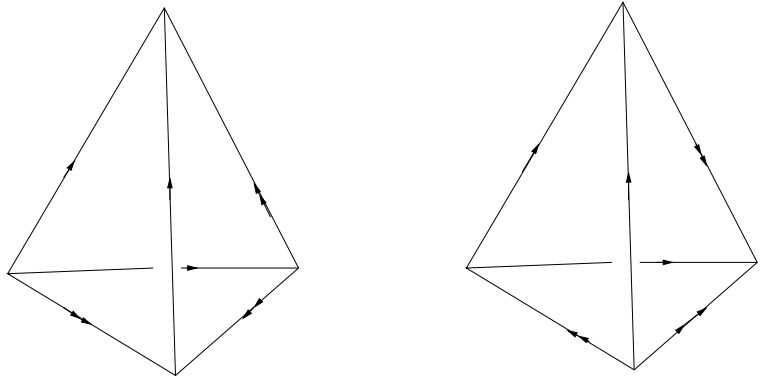}}
\end{picture}
\caption{\sl The eight figure knot}\label{figure:eight}
\end{figure}

We refer again to the parametrization of tetrahedra using
 $z_i$, $z_i'$ and $\tilde z_i$.  In this section we  obtain
 the general equations for gluing two tetrahedra according to the
scheme in Figure \ref{figure:eight}.

The first set of equations concerns the compatibility of 
Cartan's invariants of each of the four triples of points in the
tetrahedron:
$$
A\leftrightarrow A' 
\Longrightarrow\ \A(p_1,p_2,q_1)=\A(\dot p_1,\dot p_2, \dot q_2)
 \Longrightarrow t=\dot s
$$
$$
B\leftrightarrow B' 
\Longrightarrow\ \A(p_1,q_1,q_2)=\A(\dot q_1,\dot q_2, \dot p_2)
 \Longrightarrow \frac{s|z|^2-t+2\Im z}{|z-1|^2}= 
{\mbox {function of}}\  \dot z, \dot t, \dot s
$$
$$
C\leftrightarrow C' 
\Longrightarrow\ \A(p_1,p_2,q_2)=\A(\dot q_1,\dot q_2, \dot p_1)
 \Longrightarrow s=\frac{\dot s|\dot z|^2-\dot t+2\Im \dot z}{|\dot z-1|^2}
$$
$$
D\leftrightarrow D' 
\Longrightarrow\ \A(p_2,q_1,q_2)=\A(\dot p_2,\dot q_1, \dot p_1)
 \Longrightarrow {\mbox {function of}}\   z,  t,  s =\dot t
$$
The function of $z$, $t$ and $s$ above is
$$
\frac{2(s-t)\Re z +2(1+ts)\Im z +t(1+s^2)|z|^2-s(1+t^2)}{|(s-i)z+i-t|^2}.
$$
There are three independent equations, the fourth one being
 a consequence of the cocycle condition.  We choose the first 
and the last two equations.  From the first equation, $\dot s$ is
 dermined by $t$.  From the last, $\dot t$ is determined by $z$, $t$ 
and $s$.  Substituting in the third equation we obtain
$s$ as a function
of $z$, $\dot z$ and $t$.  That gives a 5 parameter family
of a couple of tetrahedra with compatible Cartan's invariants
under the gluing scheme. 

\subsection{Symmetric tetrahedra}
\begin{prop} If two symmetric tetrahedra are glued following the
scheme to obtain the complement of the figure eight knot then
they are both regular. In that case, a couple of regular tetrahedra is
parametrized by a hypersurface in the variables $z$ and $\dot z$.
\end{prop}
\Pf
For a symmetric tetrahedron $s=t$.  From the equations above we
obtain that the four triples have the same Cartan's invariant
, therefore they are regular.
In this case we have
$$
t=\frac{\Im z}{1-\Re z}=\frac{\Im \dot z}{1-\Re \dot z}
$$
\EPf

In order to have a coherent gluing of the tetrahedra along the edges
we have to impose the following equations (where we changed the dot notation 
for a different variable $w$) two for each cycle of edges,
corresponding to the two end points of each cycle:
$$
\begin{matrix}
 & z_1w_1\tilde z'_2w_3 z_2 \tilde w_1=1\\
 & z'_1 w'_1 z'_2\tilde w'_3\tilde z_2 \tilde w'_1=1\\
 & z_3\tilde w_3 \tilde z_3\tilde w'_2\tilde z_1\tilde w_2=1\\
 & \tilde z'_3 w'_3 z'_3 w'_2 \tilde z'_1 w_2=1
\end{matrix}
$$
The product of the four equations is clearly 1, so only three
of the equations are independent.
Using the relations between the invariants we simplify to
 $$
\begin{matrix}
 & {(z_2-1)\tilde z'_2(w_1-1)\tilde w_1}=1\\
 & {(z'_2-1)\tilde z_2(\tilde w'_1-1) w'_1}=1\\
 & {(\tilde z_1-1)z_3(\tilde w_3-1)\tilde w'_2}=1\\
 & {(\tilde z'_1-1) z'_3(w'_3-1) w_2}=1\\
\end{matrix}
$$
\begin{prop}
The only symmetric tetrahedra with identifications as the scheme
above giving the eight knot complement is the one obtained in 
the previous section.
\end{prop}
\Pf
If the tetrahedra are regular $z_1=\tilde z_1$.  We obtain then
$$
\tilde z'_2=\frac{1}{1-\tilde z'_1}=\frac{1}{1-\frac{1}{{\bar z}_1}}=
\frac{\bar z_1}{\bar z_1-1}=1-\bar z_2.
$$
%Therefore, the first equation becomes $|z_2-1|^2(w_1-1)w_1=-1$.
The second equation becomes ${-\bar z_1\over |1-z_1|^2}{1-\bar w_1\over \bar w_1}=1$ and the last one  ${1-\bar z_1\over \bar z_1}{-\bar w_1\over 1- w_1}=1$.
From those two equations follows that $w_1^2+\bar w_1=0$ which has the unique
solution $w_1=e^{i\pi/3}$.
\EPf
\section{Gluing special tetrahedra: The Whitehead Link}
Using the other special tetrahedra defined in \ref{sub:whitehead} we obtain
the complement of the Whitehead link.  It suffices to observe that we can
glue four tetrahedra as in Thurston forming an octahedra with
dihedral angles equal to $\pi/2$.  
We make $p_1=(0,1+\sqrt{2})$, $p_2=(0,-(1+\sqrt{2}))$, $q_1=(1,0)$
and $q_2=(i,0)$. We have  $z_1=\tilde z_1=i$ with $\A (p_1,q_1,q_2)=\pi/4$. 
We want to show completeness. Define
$q_3=(-1,0)$ and $q_4=(-i,0)$. 

\begin{figure}
\setlength{\unitlength}{1cm}
\begin{picture}(13,10)
\psfrag{Id}{Id}
\psfrag{w=q_1}{$w=q_1$}
\psfrag{w=q_3}{$w=q_3$}
\psfrag{v=q_4}{$v=q_4$}
\psfrag{v=q_2}{$v=q_2$}
\psfrag{v=p_1}{$v=p_1$}
\psfrag{v=p_2}{$v=p_2$}
\psfrag{w}{w}
\psfrag{a}{a}
\psfrag{b}{b}
\psfrag{c}{c}
\psfrag{d}{d}
\psfrag{e}{e}
\psfrag{f}{f}
\psfrag{g}{g}
\psfrag{h}{h}
\psfrag{g_1}{$G_1$}
\psfrag{g_2}{$G_2$}
\psfrag{g_3}{$G_3$}
\psfrag{g_2^{-1}}{$G_2^{-1}$}
\psfrag{g_1^{-1}}{$G_1^{-1}$}
\psfrag{g_3^{-1}}{$G_3^{-1}$}
\psfrag{tz1}{$\tilde z_1$}
\psfrag{tz3}{$\tilde z_3$}
\psfrag{tz2}{$\tilde z_2$}
\psfrag{w1}{$w_1$}
\psfrag{w2}{$w_2$}
\psfrag{w3}{$w_3$}
\psfrag{z1}{$z_1$}
\psfrag{z3}{$z_3$}
\psfrag{z2}{$z_2$}
\psfrag{tw1}{$\tilde w_1$}
\psfrag{tw3}{$\tilde w_3$}
\psfrag{tw2}{$\tilde w_2$}
\put(3,1){\includegraphics[height=8cm,width=10cm]{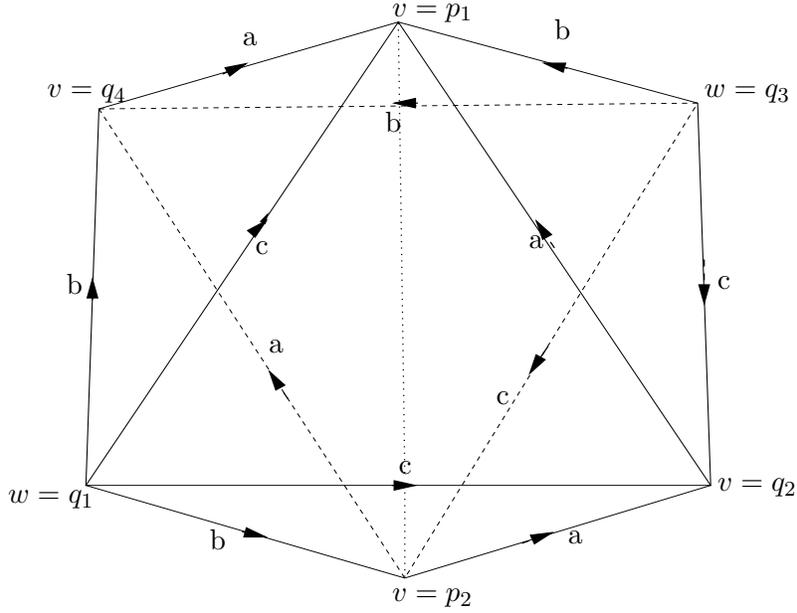}}
\end{picture}
\caption{\sl The Whitehead link complement}\label{figure:white}
\end{figure}

The generators of the group are given by
$$
g_A:[p_1,q_1,q_2]\rightarrow [q_2,q_3,p_2]
$$
$$
g_B:[p_1,q_2,q_3]\rightarrow [q_4,p_2,q_3] 
$$
$$
g_C:[p_1,q_3,q_4]\rightarrow [q_4,q_1,p_2]
$$
$$
g_D:[p_1,q_4,q_1]\rightarrow [q_2,p_2,q_1]
$$
conjugating the generators above with the mapping
$$
[p_1,q_1,q_2]\rightarrow [\infty,0,(1,1)]
$$
we obtain the following matrices in $SU(2,1)$ representing the generators:
$$
G_1=\left[\begin{matrix} 1  &  0  &  -i\\
                         -1-i  &   1      &  -1+i\\
		 	 -1-i  &   1-i &  i
\end{matrix}\right]
$$
$$
G_2=\left[\begin{matrix} 1  &  1-i  &  -1+i\\
                         -1-i  &   -1      &    1-i\\
		 	-1+i   &   1+i &  -1-2i
\end{matrix}\right]
$$
$$
G_3=\left[\begin{matrix} i        &  1+i        &  -i     \\
                         1-i      &   -1-2i     & 2i            \\
		 	-1-i      &   -3+i      & 3+2i 
\end{matrix}\right]
$$
$$
G_4=\left[\begin{matrix} -i    &  0  &  0     \\
                         -1+i  &  -1 & 0            \\
		 	-1+i   &  -1+i & -i 
\end{matrix}\right]
$$
$G_1$ and $G_3$ have trace $2+i$ and therefore are loxodromic, $G_2$ and 
$G_4$ have trace $-1-2i$ and  are elliptic of order four.

We obtained the following
\begin{thm}\label{main1}
The representation of the fundamental group of the 
Whitehead link complement generated by $G_1,G_2,G_3,G_4$ is 
in $PU(2,1,\Z[i])$ and is therefore discrete.
\end{thm}

\subsection{Holonomy}
There are two tori. We use the notation as in \cite{Ra}. We compute
 their holonomy as in the case of the figure eight knot.
 The first torus has holonomy generated by
$$
H_1=G_3^{-1} G_1^{-1}=\left[\begin{matrix} -1-6i    &  -6-4i  &  2+4i     \\
                         -4+6i  &  1+8i & 2-4i            \\
		 	2+4i   &  4+2i & -1-2i 
\end{matrix}\right]\ \ \ {\mbox and}\ \  H_2=G_2.
$$
Observe that $H_1$ is parabolic but $H_2$ is elliptic.
The other torus has holonomy generated by
$$
H'_1=G_3G_1^{-2}G_3=\left[\begin{matrix} 5    &  2-6i  &  -4     \\
                         -8-4i  &  -7+8i & 6+2i            \\
		 	-8+8i   &  8+12i & 5-8i
\end{matrix}\right]\ \ \ {\mbox and}\ \  H'_2=Id.
$$
Here $H'_1$ is parabolic.  Note that the holonomy of that torus
is not faithful.

\end{document}